\documentclass{amsart}

\newcommand{\bdis}{\begin{displaymath}}
\newcommand{\edis}{\end{displaymath}}
\newcommand{\be}{\begin{equation}}
\newcommand{\ee}{\end{equation}}

\newcommand{\mcal}{\mathcal}

\newtheorem{theorem}{Theorem}

\theoremstyle{definition}

\newtheorem{rem}[theorem]{Remark}
\newtheorem{form}[theorem]{Formula}
\newtheorem{cor}[theorem]{Corollary}

\theoremstyle{remark}

\numberwithin{equation}{section}

\begin{document}

\title{On some properties of Riemann zeta function on critical line\footnote{paper published in ACTA ARITHMETICA, XXVI (1974).}}

\author{Jan Moser}

\address{Department of Mathematical Analysis and Numerical Mathematics, Comenius University, Mlynska Dolina M105, 842 48 Bratislava, SLOVAKIA}

\email{jan.mozer@fmph.uniba.sk}


\begin{abstract}
The aim of this paper is to show further results following those published in [5], and to relate the Riemann zeta function to the
relativistic cosmology. 
\end{abstract}

\maketitle

Let $0<\gamma'<\gamma''$ be the coordinates of neighboring zeros of the zeta function:
\bdis
\rho'=\frac{1}{2}+i\gamma', \quad \rho''=\frac{1}{2}+i\gamma'' ,
\edis
and let $\{ t_0\}$ be a sequence of positive numbers such that the following three properties hold true:

\begin{equation*}
\gamma'<t_0<\gamma'', \tag{a}
\end{equation*}
\begin{equation*}
Z'(t_0)=0 \tag{b},
\end{equation*}
\begin{equation*}
t_0\to+\infty \tag{c} .
\end{equation*}
Let $\{ \tilde{t}_0\}$ be such a subsequence of the sequence $\{ t_0\}$ that
\bdis
\left|\zeta\left(\frac{1}{2}+it_0\right)\right|>\frac{1}{\tilde{t}_0^\alpha}, \quad 0<\alpha\leq 1 .
\edis
Let $\tilde{\gamma}',\ \tilde{\gamma}''$ be coordinates of such neighboring zeros of the Riemann zeta function $\zeta(s)$ that
$\tilde{\gamma}'<\tilde{t}_0<\tilde{\gamma}''$. The symbol $\{ \tilde{\gamma}',\tilde{\gamma}''\}$ will denote any sequence of such neighboring
coordinates. \\

Numerical experiments concerning the $Z(t)$ function show that the points $t_0$ are distributed with no significant dispersion in neighborhoods
of the points $(\gamma'+\gamma'')/2$. Let us denote
\bdis
\Delta(t_0)=\min \{ t_0-\gamma',\gamma''-t_0\} .
\edis
It is still possible that also in the case
\bdis
\gamma''-\gamma'>\frac{1}{\gamma^{\prime\alpha}}
\edis
$\Delta(t_0)$ will be arbitrary small, i.e. the point $t_0$ is placed arbitrarily close to $\gamma'$ or to $\gamma''$. More exactly: it arises a question
to estimate $\Delta(t_0)$ from bellow. In this the following Theorem is helpful:

\begin{theorem}
If the Riemann conjecture holds true then
\bdis
\Delta(\tilde{t}_0)>\frac{1}{\tilde{\gamma}^{\prime\alpha}}, \quad 0<\alpha\leq 1 ,
\edis
(starting with some $\tilde{t}_0$).
\end{theorem}
The proof of this Theorem can be done in the same way as the proof of the Theorem 1 in \cite{jm1}. The only difference is that the role of the
relation \cite{jm1}, (8) plays now the relation
\bdis
\Delta(\tilde{\tilde{t}}_0)\leq \frac{1}{\tilde{\tilde{\gamma}}^{\prime\alpha_0}} ,
\edis
and the relation \cite{jm1}, (5) is replaced by
\bdis
\sum_\gamma \frac{1}{\left( \tilde{\tilde{t}}_0-\gamma\right)^2}>\frac{1}{\left( \tilde{\tilde{t}}_0-\tilde{\tilde{\gamma}}^{\prime}\right)^2}\geq
\tilde{\tilde{\gamma}}^{\prime 2\alpha_0}, \quad \Delta(\tilde{\tilde{t}}_0)=\tilde{\tilde{t}}_0-\tilde{\tilde{\gamma}}^\prime ,
\edis
or
\bdis
\sum_\gamma \frac{1}{\left( \tilde{\tilde{t}}_0-\gamma\right)^2}>\frac{1}{\left( \tilde{\tilde{t}}_0-\tilde{\tilde{\gamma}}^{\prime\prime}\right)^2}\geq
\tilde{\tilde{\gamma}}^{\prime 2\alpha_0}, \quad \Delta(\tilde{\tilde{t}}_0)=\tilde{\tilde{\gamma}}^{\prime\prime}-\tilde{\tilde{t}}_0 .
\edis

\begin{rem} {\rm
The fact we are interested in the values of $\alpha$ ranged in the interval $(0,1]$ only is not important. Theorem 1 in \cite{jm1} and our Theorem 1
remain valid for $0<\alpha<+\infty$.
}
\end{rem}
We will show the following formula is true.

\begin{form} \label{form1} {\rm
If the Riemann conjecture is true then
\begin{equation*}
\sum_\gamma \frac{1}{(t-\gamma)^2}=-\frac{{\rm d}}{{\rm d}t}\left\{\frac{Z'(t)}{Z(t)}\right\}+\mcal{O}\left(\frac{1}{t}\right) . \tag{1}
\end{equation*}
}
\end{form}

This formula is followed by series of consequences concerning some features of the $Z(t)$ function, namely:

\begin{cor}
The function $Z'(t)/Z(t)$ is decreasing on the interval $\gamma'<t<\gamma^{\prime\prime}$ .
\end{cor}

\begin{cor}
The first derivative $Z'(t_0)$ and the second derivative $Z''(t_0)$ cannot vanish simultaneously.
\end{cor}

\begin{cor}
Within the interval $\gamma'<t<\gamma^{\prime\prime}$: if $Z(t)>0$ then $Z$ cannot attain two maxima and one minimum; if $Z(t)<0$ then $Z$ cannot
attain two minima and one maximum.
\end{cor}

We will give the proof of Formula \ref{form1}. It holds true:
\bdis
\zeta\left(\frac{1}{2}+it\right)=e^{-i\vartheta(t)}Z(t) ,
\edis
\bdis
\frac{\zeta'\left(\frac{1}{2}+it\right)}{\zeta\left(\frac{1}{2}+it\right)}=\vartheta'(t)-i\frac{Z'(t)}{Z(t)},
\edis
\begin{equation*}
\left\{\frac{\zeta'\left(\frac{1}{2}+it\right)}{\zeta\left(\frac{1}{2}+it\right)}\right\}^2=\vartheta^{\prime 2}(t)-
\left\{\frac{Z'(t)}{Z(t)}\right\}^2+i2\vartheta'(t)\frac{Z'(t)}{Z(t)}, \tag{2}
\end{equation*}
\bdis
\frac{Z''(t)}{Z(t)}-\left\{\frac{Z'(t)}{Z(t)}\right\}^2=i\vartheta''(t)-\frac{\zeta''\left(\frac{1}{2}+it\right)}{\zeta\left(\frac{1}{2}+it\right)}+
\left\{\frac{\zeta'\left(\frac{1}{2}+it\right)}{\zeta\left(\frac{1}{2}+it\right)}\right\}^2,
\edis
and subsequently
\begin{equation*}
\frac{\zeta''\left(\frac{1}{2}+it\right)}{\zeta\left(\frac{1}{2}+it\right)}=-\frac{Z''(t)}{Z(t)}+\vartheta^{\prime 2}(t)+i
\left[ \vartheta''(t)+2\vartheta'(t)\frac{Z'(t)}{Z(t)}\right] .
\end{equation*}

On the other hand, see \cite{jm1},
\begin{align*}
& \frac{\zeta''\left(\frac{1}{2}+it\right)}{\zeta\left(\frac{1}{2}+it\right)} & & = & & \sum_{\gamma}\frac{1}{(t-\gamma)^2}+
\left\{\frac{\zeta'\left(\frac{1}{2}+it\right)}{\zeta\left(\frac{1}{2}+it\right)}\right\}^2+\mcal{O}\left(\frac{1}{t}\right) & \tag{4} \\
& & & = & & \sum_{\gamma}\frac{1}{(t-\gamma)^2}+\vartheta^{\prime 2}(t)-\left\{\frac{Z'(t)}{Z(t)}\right\}^2+
i2\vartheta'(t)\frac{Z'(t)}{Z(t)}+\mcal{O}\left(\frac{1}{t}\right) &
\end{align*}
if used (2). Using the fact that $\vartheta''(t)=\mcal{O}(t^{-1})$, see e.g. \cite{titch}, page 260, we can compare Eqs. (3) and (4) to obtain
\bdis
\sum_{\gamma}\frac{1}{(t-\gamma)^2}=\left\{\frac{Z'(t)}{Z(t)}\right\}^2-\frac{Z''(t)}{Z(t)}+\mcal{O}\left(\frac{1}{t}\right),
\edis
that is exactly (1). \\

Eq. (4) is followed by

\begin{cor}
\bdis
\zeta''\left(\frac{1}{2}+it\right)\not=0, \qquad \gamma'<t<\gamma'' .
\edis
\end{cor}
In fact, if $\zeta''(1/2+i\tilde{t})=0,\ \gamma'<\tilde{t}<\gamma''$ then eq. (4) implies the following:
\begin{equation*}
\sum_\gamma\frac{1}{(\tilde{t}-\gamma)^2}+\vartheta^{\prime 2}(\tilde{t})-\left\{\frac{Z'(\tilde{t})}{Z(\tilde{t})}\right\}^2=
\mcal{O}\left(\frac{1}{\tilde{t}}\right) , \tag{5}
\end{equation*}
\begin{equation*}
2\vartheta'(\tilde{t})\frac{Z'(\tilde{t})}{Z(\tilde{t})}=\mcal{O}\left(\frac{1}{\tilde{t}}\right) . \tag{6}
\end{equation*}

Using the asymptotics (see e.g. \cite{titch}, p. 260)
\begin{equation*}
\vartheta'(\tilde{t})\sim \frac{1}{2}\ln(\tilde{t}) , \tag{7}
\end{equation*}
and (6) we can write
\begin{equation*}
\frac{Z'(\tilde{t})}{Z(\tilde{t})}=\mcal{O}\left(\frac{1}{\tilde{t}\ln(\tilde{t})}\right) . \tag{8}
\end{equation*}
Now, Eqs. (5) and (8) imply
\bdis
\sum_{\gamma}\frac{1}{(\tilde{t}-\gamma)^2}+\vartheta^{\prime 2}(\tilde{t})=\mcal{O}\left(\frac{1}{\tilde{t}}\right) .
\edis
This, however, is in contradiction with (7). \\

We will prove in what follows that the Riemann conjecture implies the next two {\bf asymptotic formulae} hold true:
\begin{equation*}
\sum_{\gamma}\frac{1}{(t_0-\gamma)^2} \sim  2\sum_{\gamma}\frac{t_0^2}{(t_0^2-\gamma^2)^2}-\frac{\pi}{4t_0}, \tag{A}
\end{equation*}
\begin{equation*}
\sum_{\gamma}\frac{t_0^2}{(t_0^2-\gamma^2)^2}\sim \sum_{\gamma}\frac{\gamma^2}{(t_0^2-\gamma^2)^2}-\frac{\pi}{4t_0}, \tag{B}
\end{equation*}
as $t_0\to+\infty$. \\

In the work \cite{jm2} (\cite{dav}, pages 91-92) distinguishing the real and the imaginary part of the equation:
\bdis
\frac{\zeta'\left(\frac{1}{2}+it_0\right)}{\zeta\left(\frac{1}{2}+it_0\right)}+\vartheta'(t_0)=0 ,
\edis
and using the formula:
\bdis
\frac{\zeta'(s)}{\zeta(s)}=b-\frac{1}{s-1}-\frac{{\rm d}}{{\rm d}s}\ln \Gamma\left(\frac{s}{2}+1\right)+
\sum_{\rho}\left(\frac{1}{s-\rho}+\frac{1}{\rho}\right) ,
\edis
we obtain at $t_0\to+\infty$ the Riemann formula
\bdis
\ln\left(\frac{e^{C+2}}{4\pi}\right)=\sum_{\gamma}\frac{1}{\frac{1}{4}+\gamma^2} ,
\edis
(where $C$ is the Euler constant) and on the other hand we have
\bdis
\lim_{t_0\to\infty}\sum_{\gamma}\frac{t_0}{\gamma^2-t_0^2}=\frac{\pi}{4} ,
\edis
and this this means that
\begin{equation*}
\sum_{\gamma}\frac{1}{\gamma^2-t_0^2}\sim \frac{\pi}{4t_0} . \tag{9}
\end{equation*}

Now, if we use (9) appropriately, we obtain:
\begin{align*}
& \sum_{\gamma}\frac{1}{(t_0-\gamma)^2} & & = & & \sum_{\gamma>0}\frac{1}{(t_0-\gamma)^2} + \sum_{\gamma<0}\frac{1}{(t_0-\gamma)^2}& \tag{A} \\
& & & = & & \sum_{\gamma>0}\frac{1}{(t_0-\gamma)^2}+\sum_{\gamma>0}\frac{1}{(t_0+\gamma)^2} =
\sum_{\gamma>0}\frac{(t_0-\gamma)^2+(t_0+\gamma)^2}{(t_0^2-\gamma^2)^2}& \\
& & & = & & 2\sum_{\gamma>0}\frac{t_0^2+\gamma^2}{(t_0^2-\gamma^2)^2}=2\sum_{\gamma>0}\frac{t_0^2-\gamma^2+2\gamma^2}{(t_0^2-\gamma^2)^2} &  \\
& & & = & & 2\sum_{\gamma>0}\frac{1}{t_0^2-\gamma^2}+4\sum_{\gamma>0}\frac{\gamma^2}{(t_0^2-\gamma^2)^2} & \\
& & & = & & \sum_{\gamma}\frac{1}{t_0^2-\gamma^2}+2\sum_{\gamma}\frac{\gamma^2}{(t_0^2-\gamma^2)^2}
\sim-\frac{\pi}{4t_0}+2\sum_\gamma\frac{\gamma^2}{(t_0^2-\gamma^2)^2},  &
\end{align*}

\begin{align*}
& \sum_{\gamma}\frac{t_0^2-\gamma^2}{(t_0^2-\gamma^2)^2}=\sum_{\gamma}\frac{1}{t_0^2-\gamma^2}\sim-\frac{\pi}{4t_0},  &  & & & & & & & & & & & & & & & & & & & & & & & & & & & && & & &  \tag{B}
\end{align*}
i.e.
\bdis
\sum_\gamma \frac{t_0^2}{(t_0^2-\gamma^2)^2}\sim \sum_{\gamma}\frac{\gamma^2}{(t_0^2-\gamma^2)^2}-\frac{\pi}{4t_0} .
\edis

\subsection*{On use of the Riemann zeta function in relativistic cosmology}

It is well-known (see e.g. \cite{mak}, page 209, (8.210), (8.211)) that the essence of the relativistic cosmology is in the so-called Friedmann
equations:
\begin{equation*}
\kappa c^2\rho=\frac{3}{R^2}\left(k c^2+R^{\prime 2}\right) , \tag{10}
\end{equation*}
\begin{equation*}
\kappa p = -\frac{2R''}{R}-\frac{R^{\prime 2}}{R^2}-\frac{kc^2}{R^2} , \tag{11}
\end{equation*}
(we do not consider the cosmological constant), where: $\kappa$ is the Newton gravitational constant, $c$ is the speed of light in vacuum,
$\rho(t)$ is the mass density, $p(t)$ is pressure, $k=-1,0,1$, $R(t)$ is the "radius"\ of the Universe; prime denotes the derivative with respect to
the time. \\

We have a pair of equations for a triple of unknown functions $\rho(t),p(t),R(t)$. We have to add one more equation, usually it is the so-called
\emph{state equation} of the form:
\begin{equation*}
F(\rho,p)= 0 . \tag{12}
\end{equation*}
Often used state equation reads
\begin{equation*}
p=0 , \tag{13}
\end{equation*}
that describes an universe filled with matter with vanishing pressure. \\

It is also studied an universe in which
\begin{equation*}
p=\frac{\rho c^2}{3} , \tag{14}
\end{equation*}
(see e.g. \cite{ll}, pages 387, 388). This corresponds to the case of maximal pressure $p$ at given matter density $\rho$. \\

Let us remark that for physical reasons\footnote{Let us remark this work was written before 1974.} we require
\begin{equation*}
\rho>0, \quad p\geq 0 . \tag{15}
\end{equation*}

From the mathematical point of view one constructs the models of the Universe in the following way: equation (12) is postulated and subsequently the pair
of equations (10), (11) is solved for two independent variables. \\

We will construct (supposing the Riemann conjecture holds true) an infinite set of the models of the Universe prescribing the radius $R$ as a
function of time and the conditions (15) will be fulfilled in some intervals of time. \\

We suppose:
\begin{equation*}
R(t)=|Z(t)|, \quad k=+1 , \tag{16}
\end{equation*}
this means the radius $R$ of the Universe is related with the function $\zeta(1/2+it)$ and one supposes the spherical geometry. \\

In this case using formulae (1) and (10), (11) we have

\begin{equation*}
\frac{\kappa c^2}{3}\rho(t)=\frac{c^2}{Z^2(t)}+\left\{\frac{Z'(t)}{Z(t)}\right\}^2 , \tag{17}
\end{equation*}
\begin{equation*}
\frac{\kappa}{2}p(t)=\sum_\gamma \frac{1}{(t-\gamma)^2}-\frac{3}{2}\left\{\frac{Z'(t)}{Z(t)}\right\}^2-\frac{c^2}{2}\frac{1}{Z^2(t)}+
\mcal{O}\left(\frac{1}{t}\right). \tag{18}
\end{equation*}

We will find the intervals of time $t$ in which the conditions (15) are fulfilled. \\

The first one of them - $\rho(t)>0$ - holds true for every $t>0$. \\

We will find the intervals in which also $P(t)\geq 0$ holds true. \\

From one side, using the Littlewood bound, see \cite{titch}, page 223,
\bdis
\gamma''-\gamma'<\frac{A}{\ln\ln\ln \gamma'} ,
\edis
we have
\begin{equation*}
\sum_\gamma \frac{1}{(t-\gamma)^2}>\frac{1}{(\gamma''-\gamma')^2}>A_1(\ln\ln\ln t_0)^2 . \tag{19}
\end{equation*}
From the other side, following the Littlewood - Titchmarch $\Omega$-theorem  we have that there is such a subsequence $\{\tilde{t}_0\}$ of the
sequence $\{ t_0\}$ that
\begin{equation*}
|Z(\tilde{t}_0)|>A_2\exp(\ln^\beta \tilde{t}_0), \quad \beta<\frac{1}{2} . \tag{20}
\end{equation*}

Let us remark that
\begin{equation*}
Z'(\tilde{t}_0)=0 . \tag{21}
\end{equation*}
Using (19),(20) and (21) we have there is a sequence $\{ \delta(\tilde{t}_0)\}$, $\delta(\tilde{t}_0)>0$, such that:
\begin{equation*}
\tilde{\gamma}'<\tilde{t}_0-\delta(\tilde{t}_0)<\tilde{t}_0+\delta(\tilde{t}_0)<\tilde{\gamma}^{\prime\prime} , \tag{22}
\end{equation*}
\begin{equation*}
p(t)\geq 0, \quad t\in [\tilde{t}_0-\delta(\tilde{t}_0),\tilde{t}_0+\delta(\tilde{t}_0)],
\end{equation*}
hold true. Finally, this defines the above mentioned infinite set of the models of the Universe.

\end{document}